%% file: mainIFACWC2023.tex
\documentclass{ifacconf}

\usepackage{graphicx}      
\usepackage{natbib}        
\usepackage{amsmath}
\usepackage{amssymb}
\usepackage{subcaption}
\usepackage{blkarray, bigstrut} %
\usepackage{mathtools}

\newcommand{\tabitem}{~~\llap{\textbullet}~~}

\begin{document}
\begin{frontmatter}

\title{Dynamic Optimization for \\ Monoclonal Antibody Production} 

\thanks[footnoteinfo]{Corresponding author: J.B. J{\o}rgensen (E-mail: jbjo@dtu.dk).}

\author{Morten Wahlgreen Kaysfeld}, 
\author{Deepak Kumar},
\author{Marcus Krogh Nielsen},
\author{John Bagterp Jørgensen}


\address{Department of Applied Mathematics and Computer Science,\\ 	Technical University of Denmark, DK-2800 Kgs. Lyngby, Denmark}



\begin{abstract}                

\input{tex/abstract}
\end{abstract}


\end{frontmatter}

\input{tex/introduction}

\input{tex/model}
\input{tex/optimization}

\input{tex/results}

\input{tex/conclusion}


\bibliography{bib/ifacconf,ref/References}             	
                                                   

\end{document}

%% file: tex/abstract.tex
This paper presents a dynamic optimization numerical case study for Monoclonal Antibody (mAb) production. The fermentation is conducted in a continuous perfusion reactor. We represent the existing model in terms of a general modeling methodology well-suited for simulation and optimization. The model consists of six ordinary differential equations (ODEs) for the non-constant volume and the five components in the reactor. We extend the model with a glucose inhibition term to make the model feasible for optimization case studies. We formulate an optimization problem in terms of an optimal control problem (OCP) and consider four different setups for optimization. Compared to the base case, the optimal operation of the perfusion reactor increases the mAb yield with $44$\% when samples are taken from the reactor and with $52$\% without sampling. Additionally, our results show that multiple optimal feeding trajectories exist and that full glucose utilization can be forced without loss of mAb formation.





%% file: tex/introduction.tex
\section{Introduction}

Monoclonal antibody (mAb) production in mammalian cells is a well-known technique for synthesising identical antibodies. These antibodies are proteins with significance in medical applications. In 2017, mAbs represented the 6 top-selling biopharmaceutical products and has an expected yearly sale of 130-200 billion US dollars in 2022 \citep{walsh:2018a,grilo:2019a}. This has resulted in significant efforts to increase the synthesis of mAbs in bioreactors with mammalian cell cultures. Biopharmaceutical companies increasingly seek for novel protein production methods to accommodate for an increasing number of protein drug candidates that enter various phases of research. In today's competitive market, it is challenging to retain desirable quality attributes while shortening time to market, maintaining cost efficiency, and enabling manufacturing flexibility.

High cell density and productivity at large scale are key factors to achieve high process yield. There exist multiple bioreactor designs to achieve a balance between cell growth and volumetric productivity \citep{blunt:2018a,mitra:2022a,carvalho:2017a}. Step-wise bolus injections of the feed solution to the production bioreactor is the most frequently applied method due to its simplicity \citep{maria:2021a}. Recently, continuous bioreactor systems have received increasing attention due to its easy scale-up, waste minimization, and non-sterile cultivation \citep{blunt:2018a}.

Medium development for a bioreactor system consists of multiple parts including optimization of feeding strategies, and production of both batch medium and feed concentrates. Thus, medium development with statistical design of experiments type methods requires lots of time, effort, and cost \citep{wohlenberg:2022a,luna:2014a,rendon:2021a}. Therefore, process models can lead to valuable insights and optimization of bioreactors without the need of performing expensive and time consuming experiments. As an example, mechanistic modeling and optimization has been applied for a U-loop reactor for single-cell protein production \citep{ritschel:2019a}. 

Mechanistic models are extensively implemented for bioreactor optimization, since the models provide a deeper understanding of the growth and production of mammalian cells \citep{glen:2018a,sha:2018a}. Glucose concentration, lactate concentration, and temperature of the cellular environment impact their metabolic pathway \citep{sissolak:2019a, fan:2015a}. A higher glucose concentration results in decreased cellular growth rate and increased product formation \citep{vergara:2018a}, whereas a higher lactate concentration deteriorates both cell growth and productivity \citep{li:2012a}.

In this paper, we consider an existing mechanistic model for mAb production and apply advanced optimization techniques to compute novel optimal feeding strategies for the process. The model describes a fermentation process for mAb production conducted in a continuous perfusion bioreactor \citep{kumar:2022a}. We present the model with a general modeling methodology well-suited for simulation and optimization \citep{wahlgreen:2022a}, and extend the model with a glucose inhibition term to ensure that optimization does not exceed physical glucose inhibition limits. The model consists of six ordinary differential equations (ODEs) for the non-constant volume and five components in the bioreactor. We present an optimization problem expressed as an optimal control problem (OCP). The OCP maximizes the final amount of mAb in the reactor while satisfying a set of operational constraints. Our results show that optimal process operation improves the mAb productivity with up to $52\%$. Additionally, our results show that there exists multiple feeding trajectories resulting in the same mAb production. Similar results have been obtained for fed-batch fermentation with Haldane growth kinetics \citep{ryde:2021a}.


The remaining part of the paper is organized as follows. Section \ref{sec:mAb} presents the model for mAb production with a general modeling methodology. Section \ref{sec:optimization} introduces the considered optimization problem formulated as an OCP. Section \ref{sec:results} presents our optimization and simulation results. Section \ref{sec:conclusion} presents our conclusions.

%% file: tex/model.tex
\section{Monoclonal Antibody Production}
\label{sec:mAb}
We consider a biotechnological process for mAb production \citep{kumar:2022a} and reformulate the model in terms of a general modeling methodology for chemical reacting systems \citep{wahlgreen:2022a}.

\subsection{General model for a continuous perfusion reactor}
The biotechnological process is conducted in a continuous perfusion reactor. We apply a general ODE model,
\begin{subequations} \label{eq:general_model}
	\begin{align} 
		\frac{dV}{dt} &= e^{\top}F_{\mathrm{in}} - F_{\mathrm{out}} - F_{\mathrm{per}}, \\
		\frac{dm}{dt} &= C_{\mathrm{in}}F_{\mathrm{in}} - c F_{\mathrm{out}} - C_{\mathrm{per}} c F_{\mathrm{per}} + RV,
	\end{align}
\end{subequations}
where $V$ is the volume, $F_{\mathrm{in}}$ is a vector with inlet flow rates, $e$ is a vector of ones of proper dimension, $F_{\mathrm{out}}$ is a scalar outlet flow rate, $F_{\mathrm{per}}$ is a scalar perfusion flow rate, $m$ is a vector of masses for each component, $C_{\mathrm{in}}$ is a matrix with inlet concentrations, $c = m/v$ is a vector with concentrations (densities), $C_{\mathrm{per}}$ is a diagonal matrix with elements between $0$ and $1$ describing the percentage of each component removed from the bioreactor by the perfusion stream, and $R$ is a vector with production rates,
\begin{align}
    R = S^{\top} r(c),
\end{align}
with $S$ being the stoichiometric matrix and $r(c)$ being a vector with reaction rates.

\subsection{Reaction stoichiometry and kinetics}
The process consists of five components,
\begin{align}
    \mathcal{C} = \{ X_v, X_d, G, L, P \},
\end{align}
where $X_v$ are viable cells, $X_d$ are dead cells, $G$ is glucose, $L$ is lactate, and $P$ is mAb (product). We express the process in terms of six stoichiometric reactions,
\begin{alignat*}{4}
	&1. \text{ Cell division},         & \alpha_{1,G}   G + X_v &\longrightarrow 2X_v + \alpha_{1,P} P,    &\hspace{0.1cm}& r_1,     \\
	&2. \text{ Cell death},            &                    X_v &\longrightarrow  X_d                 ,    && r_2,     \\ 
	&3. \text{ Maintenance 1},         & \alpha_{3,G}   G + X_v &\longrightarrow  X_v + \alpha_{3,P} P,    && r_3,     \\
    &4. \text{ Maintenance 2},         &                    X_v &\longrightarrow  X_v + \alpha_{4,L} L,    && r_4,     \\ 
    &5. \text{ Lactate production 1},  &                    X_v &\longrightarrow  X_v + \alpha_{5,L} L,    && r_5,     \\
    &6. \text{ Lactate production 2},  &                    X_v &\longrightarrow  X_v + \alpha_{6,L} L,    && r_6,
\end{alignat*}
which are compactly written in the stoichiometric matrix,
\begin{equation}
	S =
	\begin{blockarray}{*{5}{c} l}
		\begin{block}{*{5}{>{$\footnotesize}c<{$}} l}
			X$_v$    & X$_d$         & G             & L                 & P                             &                       \\
		\end{block}
		\begin{block}{[*{5}{c}]>{$\footnotesize}l<{$}}
			1        & 0             & -\alpha_{1,G} & 0                 & \alpha_{1,P}  \bigstrut[t]    & $1$ \rule{0pt}{3.4ex} \\
			-1       & 1             & 0             & 0                 & 0                             & $2$                   \\
			0        & 0             & -\alpha_{3,G} & 0                 & \alpha_{3,P}                  & $3$                   \\
			0        & 0             & 0             & \alpha_{4,L}      & 0                             & $4$                   \\
			0        & 0             & 0             & \alpha_{5,L}      & 0                             & $5$                   \\
            0        & 0             & 0             & \alpha_{6,L}      & 0             \bigstrut[b]    & $6$                   \\
		\end{block}
	\end{blockarray}
\end{equation}
The six reaction rates are given as
\begin{subequations}
    \begin{align}
        r_1 &= \mu_X(c,T)       c_{X_v}, &
        r_2 &= \mu_D(T)         c_{X_v}, \\
        r_3 &= \mu_{m_1}        c_{X_v}, &
        r_4 &= \mu_{m_2}        c_{X_v}, \\
        r_5 &= \mu_{L,p_1}(c,T) c_{X_v}, &
        r_6 &= \mu_{L,p_2}(c)   c_{X_v}, 
    \end{align}
\end{subequations}
where the different rate functions are
\begin{subequations}
    \begin{align}
        \mu_X       &= \mu_{X,max} f_{lim} f_{inh} f_{temp}             ,       \\
        \mu_D       &= \mu_{D,max} f_{D,temp}                           ,       \\
        \mu_{m_1}   &= \bar\mu_{m_1}                                    ,       \\
        \mu_{m_2}   &= \bar\mu_{m_2}\frac{L_{max,2}-c_L}{L_{max,2}}     ,       \\
        \mu_{L,p_1} &= \mu_X \frac{L_{max,1} - c_L}{L_{max,1}}          ,       \\
        \mu_{L,p_2} &= \bar \mu_{L,p_2} \frac{L_{max,1}-c_L}{L_{max,1}} ,       
    \end{align}
\end{subequations}
and
\begin{subequations} \label{eq:f_terms}
    \begin{align}
        f_{lim}     &= \frac{c_G}{K_G c_{X_v} + c_G}                        ,   \\
        f_{inh}     &= \frac{KI_L}{KI_L + c_L} (1 - KI_{P} c_P )         ,   \\
        f_{temp}    &= \exp\left( -\frac{K_1}{T} \right)                    ,   \\
        f_{D,temp}  &= \exp\left( -\frac{K_2}{T} \right)                    .
    \end{align}
\end{subequations}
The model, (\ref{eq:general_model})-(\ref{eq:f_terms}), is identical to the model presented by \cite{kumar:2022a}.


\subsection{ Model extension for optimization }
\begin{figure}[tb]
    \centering
    \includegraphics[width=0.48\textwidth]{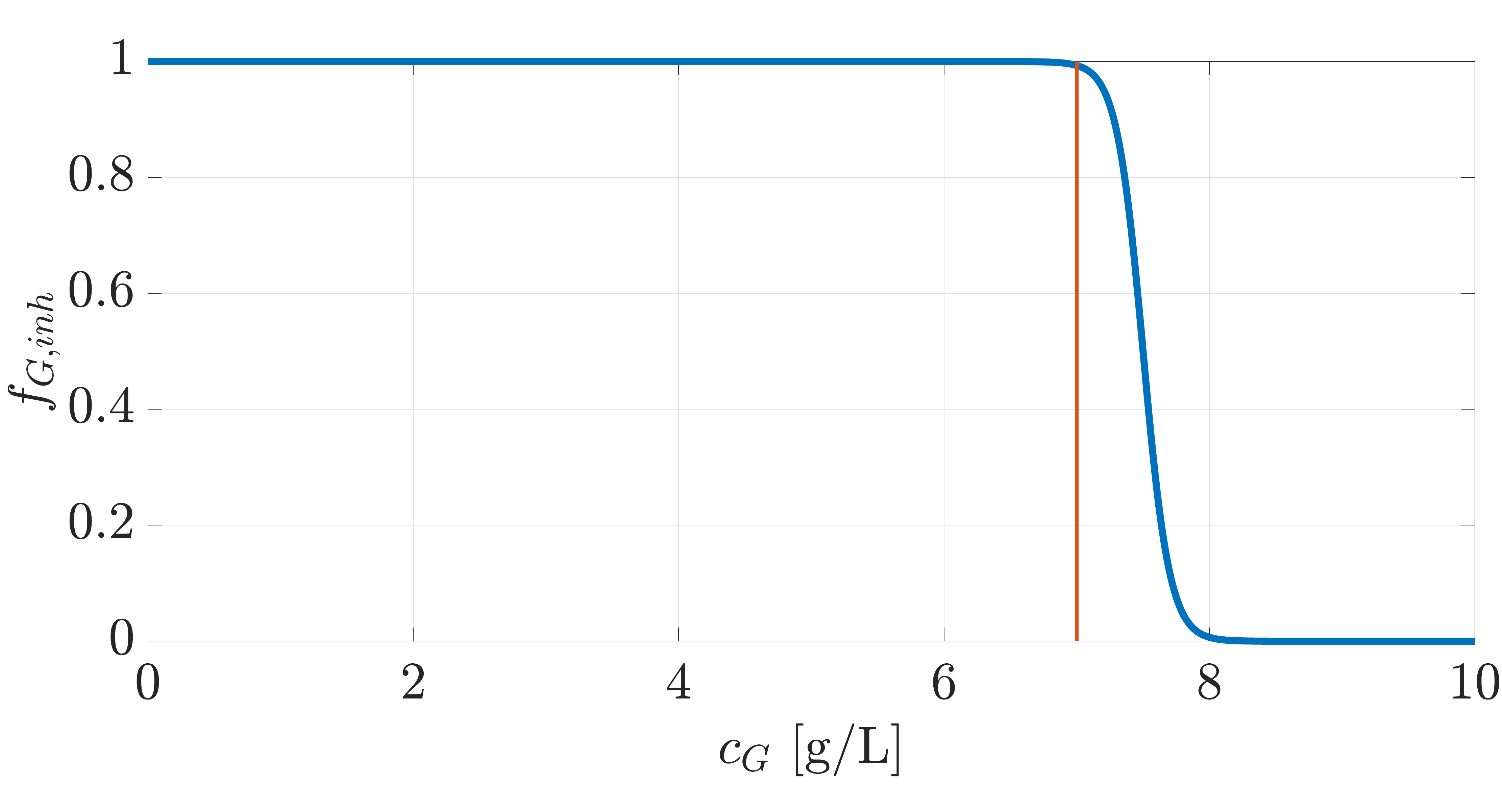}
    \caption{Glucose inhibition term, $f_{G,inh}$. Red line: $7.0$ g/L. }
    \label{fig:sigmoid}
\end{figure}

\begin{figure}[tb]
    \centering
    \includegraphics[width=0.48\textwidth]{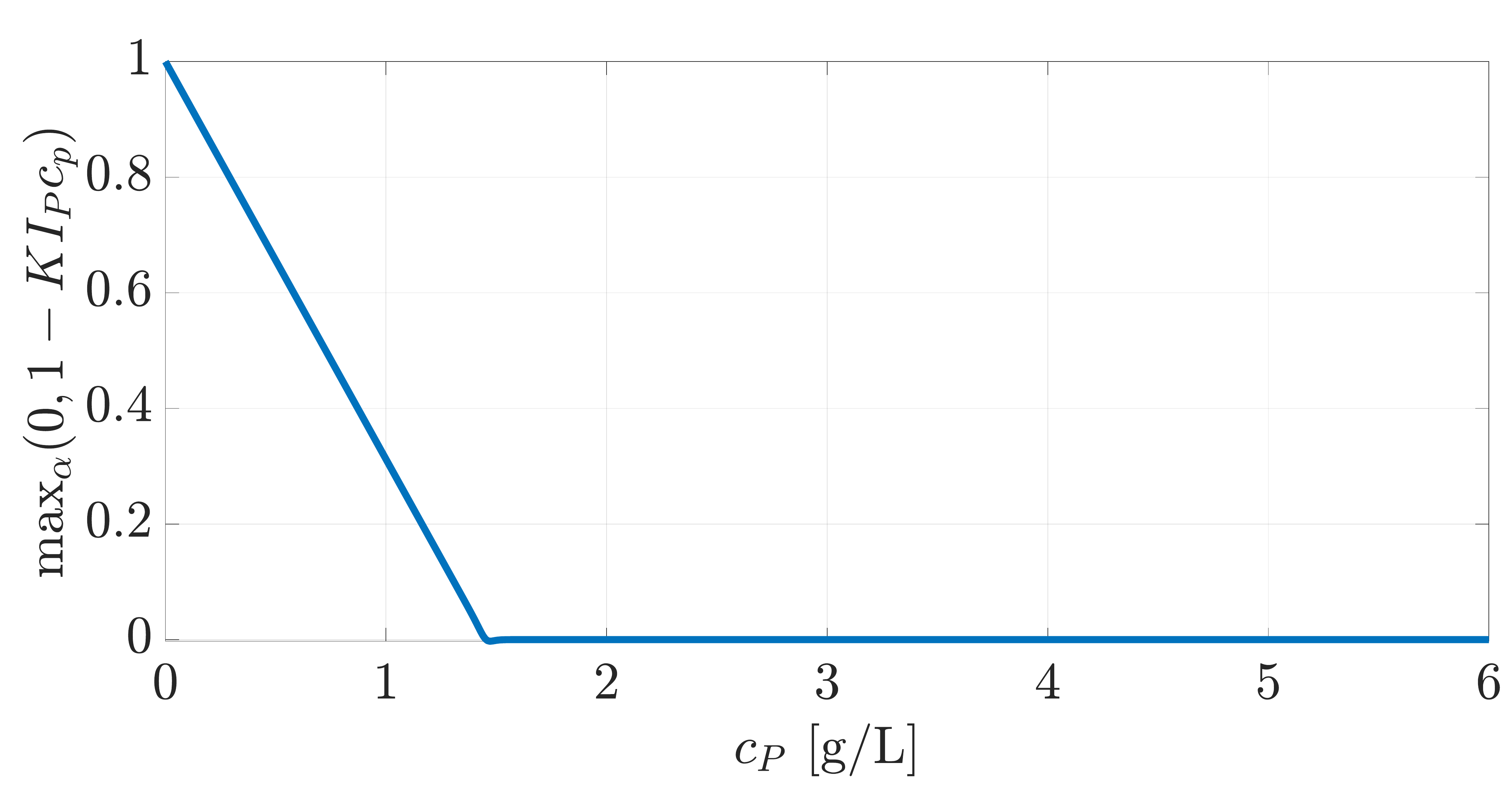}
    \caption{Smooth maximum function for product inhibition.}
    \label{fig:smax}
\end{figure}

We point out that the model, (\ref{eq:general_model})-(\ref{eq:f_terms}), does not include glucose inhibition for cell growth. In addition, the inhibition term, $f_{inh}$, becomes negative for high product concentrations. As such, the model is not directly suitable for optimization purposes. However, small model extensions lead to a feasible model well-suited for optimization. We extend the model with a simple glucose inhibition term,
\begin{align}
    f_{G,inh} &= 1 - s_{\gamma}(c_G, \bar c_{G}),
\end{align}
where $s_{\gamma}$ is a Sigmoid function given as
\begin{align}
    s_{\gamma}(x, \bar x) &= \frac{1}{1 + \exp(-\gamma(x - \bar x))}.
\end{align}
Note that $f_{G,inh} = 1$ for small glucose concentrations and rapidly converges to $0$ when $c_G > \bar c_{G}$. In addition, we remove the possibility of negative growth rate by adding a smooth maximum approximation to the product inhibition term. By combining the smooth maximum approximation with the glucose inhibition term, we get the final inhibition term,
\begin{align}
    f_{inh} &= \frac{KI_L}{KI_L + c_L} \text{max}_{\alpha}(0, 1 - KI_{P} c_P ) f_{G,inh},          
\end{align}
where $\text{max}_{\alpha}$ is a smooth maximum approximation given as
\begin{align}
    \text{max}_{\alpha}(x_1, ..., x_n) &= \frac{\sum_{i=1}^n x_i \exp(\alpha x_i)}{\sum_{i=1}^n \exp(\alpha x_i) }.
\end{align}
We design our glucose inhibition term such that the growth rate decreases for $c_G > 7$ [g/L] \citep{vergara:2018a}. As such, we select $\bar c_{G} = 7.5$ [g/L] together with $\gamma = 10.0$. Fig. \ref{fig:sigmoid} shows the activation of the glucose inhibition term, $f_{G,inh}$. Fig. \ref{fig:smax} presents the smooth maximum approximation of the product inhibition term for $\alpha = 100.0$. We point out that alternative glucose inhibition terms can be applied without loss of generality.

\subsection{ Operation of the continuous perfusion reactor }
We operate the reactor in continuous perfusion mode. We feed glucose through an inlet stream with flow rate, $F_G$, and glucose concentration, $c_{G,\mathrm{in}}$. In addition, we apply a pure water inlet stream with flow rate, $F_W$, resulting in the inlet flow vector,
\begin{align}
    F_{\mathrm{in}} = \begin{bmatrix}
        F_W \\
        F_G
    \end{bmatrix}.
\end{align}
As such, the inlet concentration matrix becomes
\begin{align}
    C_{\mathrm{in}} = \begin{bmatrix}
        0 & 0 \\
        0 & 0 \\
        0 & c_{G,\mathrm{in}} \\
        0 & 0 \\
        0 & 0 
    \end{bmatrix}.
\end{align}
The perfusion outlet has a filter that only lets spend media pass, i.e., glucose and lactate. As such, the perfusion matrix becomes,
\begin{align}
    C_{\mathrm{per}} = \begin{bmatrix}
        0 \\
        & 0 \\
        & & 1 \\
        & & & 1 \\
        & & & & 0
    \end{bmatrix}.
\end{align}
We apply the outlet stream for sampling. 



\subsection{General notation}
We formulate the model, (\ref{eq:general_model}), in terms of the general ODE system,
\begin{align} \label{eq:ODE_notation}
    \dot x(t) &= f(t, x(t), u(t), p),   & x(t_0) &= x_0,
\end{align}
where $p$ are the parameters, and the states, $x$, and the inputs, $u$, are given as
\begin{align}
    x &= \begin{bmatrix}
        V \\
        m
    \end{bmatrix}, & u &= \begin{bmatrix}
        F_{\mathrm{in}}  \\
        F_{\mathrm{per}} \\
        F_{\mathrm{out}} \\
        T
    \end{bmatrix}.
\end{align}
We apply the general ODE notation, (\ref{eq:ODE_notation}), to formulate optimization problems.


%% file: tex/optimization.tex
\section{Optimization}
\label{sec:optimization}

We formulate an optimization problem in terms of an OCP. The solution is the state and input trajectories in the finite horizon, $T_h$. We denote the initial time $t_0$ and the final time $t_f = t_0 + T_h$. We split the prediction and control horizon, $T_h$, into $N$ control intervals of equal size $T_s$. As such, $T_h = NT_s$. We assume zero order hold parameterization of the inputs,
\begin{align}
    u(t) &= u_k, && k = 0,...,N-1,
\end{align}
and consider the following OCP formulation,
\begin{subequations}
	\label{eq:OPC_general}
	\begin{alignat}{3}
		& \min_{ \xi } \quad && \varphi,                                              \\
		& \hspace{0.2cm} \text{s.t.} && x(t_0) = x_0,                                 \\
        & && \dot x(t) = f(t,x,u,d,p), \hspace{0.4cm}&& t_0 \leq t \leq t_{0}+T_h,    \\
        & && x_{\min} \leq x_{k} \leq x_{\max}, && k = 1,...,N,                       \\
		& && u_{\min} \leq u_{k} \leq u_{\max}, && k = 0,...,N-1, 
	\end{alignat}
\end{subequations} 
where $\xi = \{ u_k, x_{k+1} \}_{k=0}^{N-1}$ are the decision variables and $\varphi = \varphi(\xi)$ is the objective function. We intent to maximize the final mAb production while satisfying a set of operational constraints. This can be formulated in terms of the following OCP,
\begin{subequations}
	\label{eq:OPC}
	\begin{alignat}{3}
		& \min_{ \xi } \quad && \varphi = -m_P(t_f),                                                                          \\
		& \hspace{0.2cm} \text{s.t.} && x(t_0) = x_0,                                                                         \\
        & && \dot x(t) = f(t,x,u,d,p), \hspace{0.4cm}&& t_0 \leq t \leq t_{0}+T_h,                                            \\
        & && V_{\min} \leq V_{k} \leq V_{\max},                 && k = 1,...,N,                                               \\
        & && m_{G,\min} \leq m_{G,k}                            && k = 1,...,N,                                               \\
        & && m_{L,\min} \leq m_{L,k}                            && k = 1,...,N,                                               \\
		& && F_{\min} \leq F_{W,k} \leq F_{\max},               && k = 0,...,N-1,                                             \\
        & && F_{\min} \leq F_{G,k} \leq F_{\max},               && k = 0,...,N-1,                                             \\
        & && F_{\min} \leq F_{\mathrm{per},k} \leq F_{\max},    && k = 0,...,N-1,                                             \\
        & && T_{\min} \leq T_{k} \leq T_{\max},                 && k = 0,...,N-1.                                          
	\end{alignat}
\end{subequations}  
We point out that additional non-negativity constraints could be applied, but have not been required for meaningful optimization.

%% file: tex/results.tex
\section{Results}
\label{sec:results}

This section presents our results. We perform a base case simulation that reproduces the results from \cite{kumar:2022a}. Additionally, we perform four different optimizations and compare these to the base case. In all simulations, we consider a $14$ days fermentation.



\subsection{ Base case simulation }
\begin{figure*}
        \centering
        \includegraphics[width=1.0\textwidth]{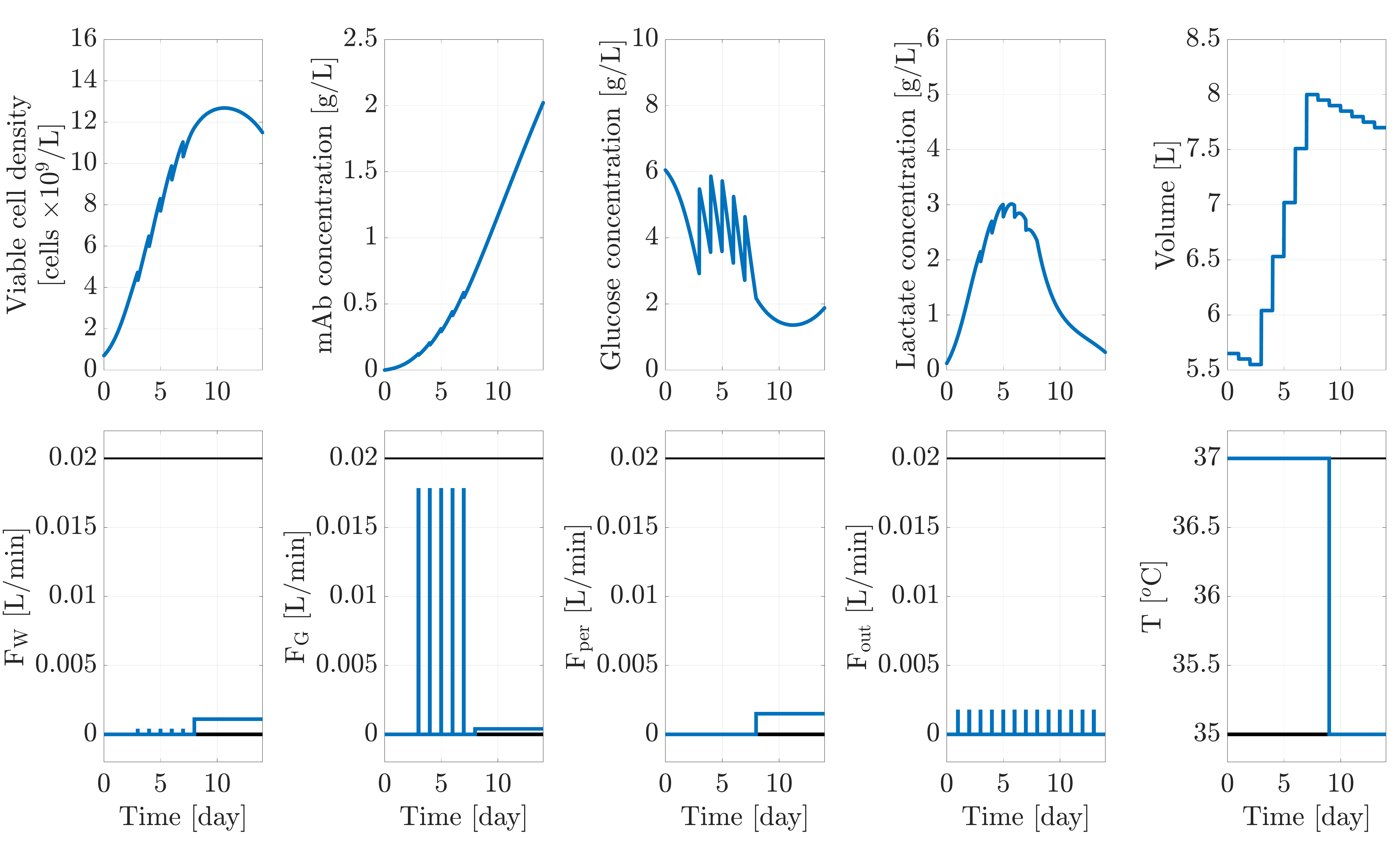}
        \caption{ Base case simulation reproducing the results from \cite{kumar:2022a}. }
        \label{fig:openLoopSim}
\end{figure*}

We simulate the model to reproduce the results presented by \cite{kumar:2022a}. As such, we operate the reactor in three phases, 1) batch phase, 2) fed-batch phase, and 3) perfusion phase. In the batch phase, no inlets or outlets are active. In the fed-batch phase, both water and glucose inlet streams are active in boluses once a day in $30$ min intervals to achieve a total inlet flow rate $F_{\mathrm{tot}} = F_W + F_G = 0.018$ L/min with glucose concentration $32.0$ g/L. In particular, this requires $F_W = 2.7692 \cdot 10^{-4}$ L/min and $F_G = 0.0177$ L/min. In the perfusion phase, the perfusion outlet is active with $F_{\mathrm{per}} = 0.0015$ L/min and the inlets are active at $F_W = 0.0011$ L/min and $F_G = 3.9923 \cdot 10^{-4}$ L/min to achieve a total inlet flow rate of $F_{\mathrm{tot}} = 0.0015$ L/min at a glucose concentration of $8.65$ g/L. Table \ref{tab:parameters} presents the list of model parameters adapted from \citet{kumar:2022a} and Table \ref{tab:operationParameters} presents the operation parameters including initial conditions. Fig. \ref{fig:openLoopSim} presents the base case simulation reproducing the results by \cite{kumar:2022a}.

\begin{table}[tb]
    \centering
    \caption{ Model parameters.}
    \begin{tabular}{lrl}
        \hline
        Parameter           & Value                         & Unit                                      \\ \hline
        $\mu_{X,max}$       & $0.153$                       & [min$^{-1}$]              \bigstrut[t]    \\
        $\mu_{D,max}$       & $3.955 \cdot 10^{-5}$         & [min$^{-1}$]                              \\
        $\bar\mu_{m_1}$     & $1.0$                         & [min$^{-1}$]                              \\
        $\bar\mu_{m_2}$     & $1.0$                         & [min$^{-1}$]                              \\
        $\bar\mu_{L,p_2}$   & $1.0$                         & [min$^{-1}$]                              \\
        $K_1$               & $1689$                        & [K]                                       \\
        $K_2$               & $524$                         & [K]                                       \\
        $K_G$               & $0.85$                        & [g/(cells$\times 10^9$)]                  \\
        $KI_L$              & $344$                         & [g/L]                                     \\
        $KI_P$              & $6.88\times 10^{-1}$          & [L/g]                                     \\
        L$_{\max,1}$        & $628$                         & [g/L]                                     \\
        L$_{\max,2}$        & $0.5$                         & [g/L]                                     \\
        $\bar c_{G}$        & $7.5$                         & [g/L]                                     \\
        $\alpha_{1,G}$      & $0.4876$                      & [g/(cells$\times 10^9$)]                  \\
        $\alpha_{1,P}$      & $6.62\times 10^{-8} $         & [g/(cells$\times 10^9$)]                  \\
        $\alpha_{3,G}$      & $1.102\times 10^{-4}$         & [g/(cells$\times 10^9$)]                  \\
        $\alpha_{3,P}$      & $1.2\times 10^{-5}$         & [g/(cells$\times 10^9$)]                  \\
        $\alpha_{4,L}$      & $1.89\times 10^{-5}$          & [g/(cells$\times 10^9$)]                  \\
        $\alpha_{5,L}$      & $0.5504$                      & [g/(cells$\times 10^9$)]                  \\
        $\alpha_{6,L}$      & $1.0249\times 10^{-5}$        & [g/(cells$\times 10^9$)]  \bigstrut[b]    \\ 
        \hline
    \end{tabular}
    \label{tab:parameters}
\end{table}

\begin{table}[tb]
    \centering
    \caption{ Operation parameters. }
    \begin{tabular}{lrl}
        \hline
        Parameter               & Value     & Unit                                  \\ \hline 
        $V_0$                   & $5.650$   & [L]                   \bigstrut[t]    \\
        $m_{X_v,0}$             & $3.955$   & [cells$\times 10^9$]                  \\
        $m_{X_d,0}$             & $0.0$     & [cells$\times 10^9$]                  \\
        $m_{G,0}$               & $34.18$   & [g]                                   \\
        $m_{L,0}$               & $0.678$   & [g]                                   \\
        $m_{P,0}$               & $0.0$     & [g]                                   \\
        $c_{G,\mathrm{in}}$     & $32.5$    & [g/L]                                 \\
        $F_{\min}$              & $0.0$     & [L/min]                               \\
        $F_{\max}$              & $0.02$    & [L/min]                               \\
        $T_{\min}$              & $308.15$  & [K]                                   \\
        $T_{\max}$              & $310.15$  & [K]                                   \\
        $V_{\min}$              & $4.0$     & [L]                                   \\
        $V_{\max}$              & $8.0$     & [L]                                   \\
        $m_{G,\min}$            & $0.0$     & [g]                                   \\
        $m_{L,\min}$            & $0.0$     & [g]                   \bigstrut[b]    \\
        \hline
    \end{tabular}
    \label{tab:operationParameters}
\end{table}


\begin{table}[tb]
    \centering
    \caption{Optimization setups. }
    \label{tab:optimizationSetups}
    \begin{tabular}{l| l l l}
        Setup   & Description                           & Constraints in OCP                                    &             \bigstrut[t]  \\ \hline
        (1)     & $\times$ sampling                     & \tabitem $F_{\mathrm{out}} = 0.0$                     & [L/min]                   \\
                & $\times$ glucose utilization          & \tabitem $m_G(t_f) \leq \infty$                       & [g]         \bigstrut[b]  \\
        (2)     & $\times$ sampling                     & \tabitem $F_{\mathrm{out}} = 0.0$                     & [L/min]                   \\
                & $\checkmark$ glucose utilization      & \tabitem $m_G(t_f) \leq 1.0$                          & [g]         \bigstrut[b]  \\
        (3)     & $\checkmark$ sampling                 & \tabitem $F_{\mathrm{out}} = \bar F_{\mathrm{out}}$   & [L/min]                   \\
                & $\times$ glucose utilization          & \tabitem $m_G(t_f) \leq \infty$                       & [g]         \bigstrut[b]  \\
        (4)     & $\checkmark$ sampling                 & \tabitem $F_{\mathrm{out}} = \bar F_{\mathrm{out}}$   & [L/min]                   \\
                & $\checkmark$ glucose utilization      & \tabitem $m_G(t_f) \leq 1.0$                          & [g]         \bigstrut[b]  \\
        \hline
    \end{tabular}   
\end{table}

\begin{table}[tb]
    \centering
    \caption{ mAb production and improvement relative to the base case. }
    \label{tab:mAbProduction}
    \begin{tabular}{lrr}
        \hline                                     
         Simulation & mAb [g]   & Improvement [\%]                  \\ \hline 
         Base case  & $15.57$   & -                 \bigstrut[t]    \\
         Opt. (1)   & $23.63$   & $52$                              \\
         Opt. (2)   & $23.63$   & $52$                              \\
         Opt. (3)   & $22.47$   & $44$                              \\
         Opt. (4)   & $22.47$   & $44$              \bigstrut[b]    \\ \hline 
    \end{tabular}
\end{table}


\begin{figure*}
        \centering
        \includegraphics[width=1.0\textwidth]{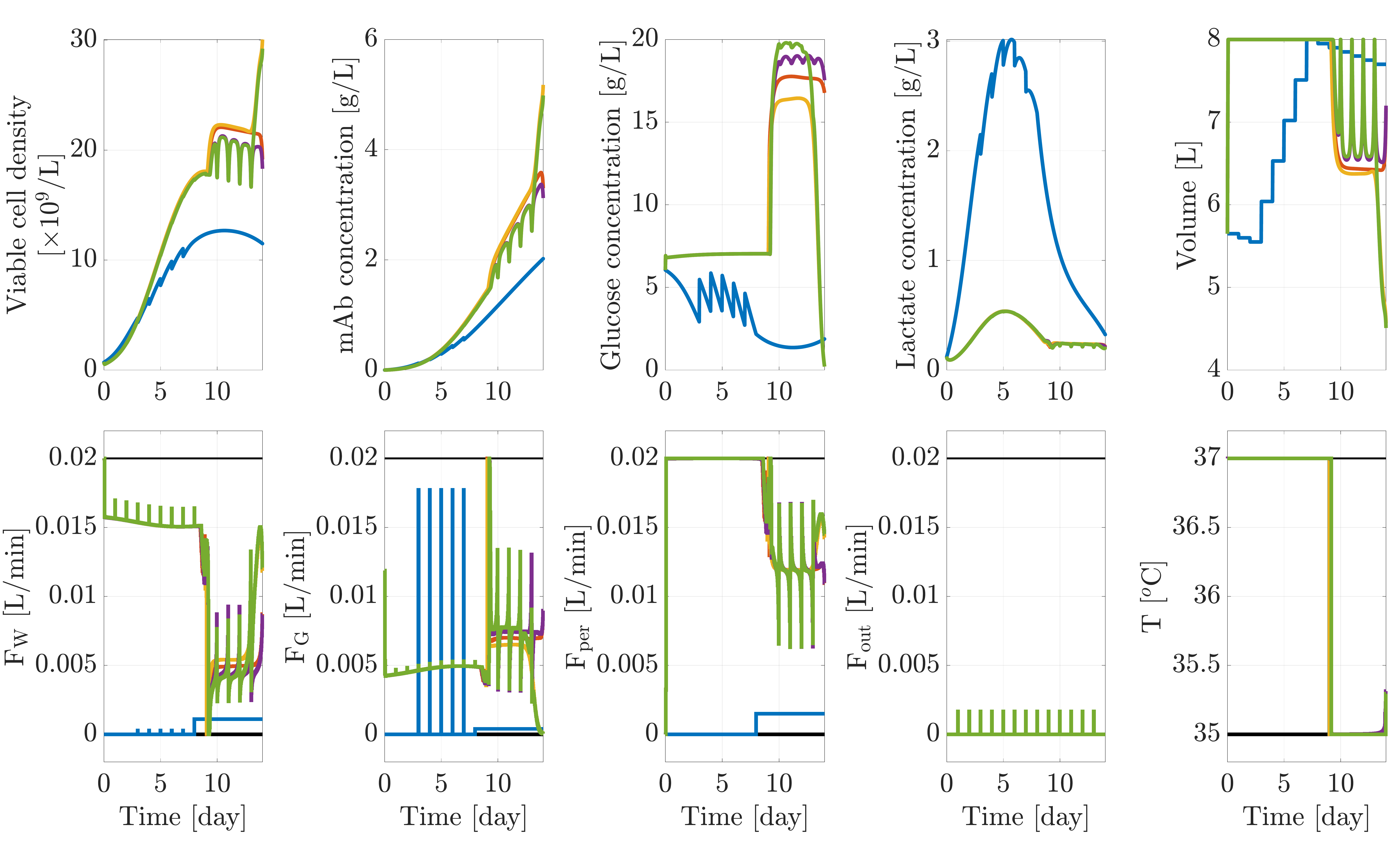}
        \caption{ The base case simulation and simulations for four different optimization setups. The blue line is the base case, the red line is setup (1), the yellow line is setup (2), the purple line is setup (3), and the green line is setup (4). }
        \label{fig:optimalSolution}
\end{figure*}

\begin{figure}[tb]
    \centering
    \includegraphics[width=0.48\textwidth]{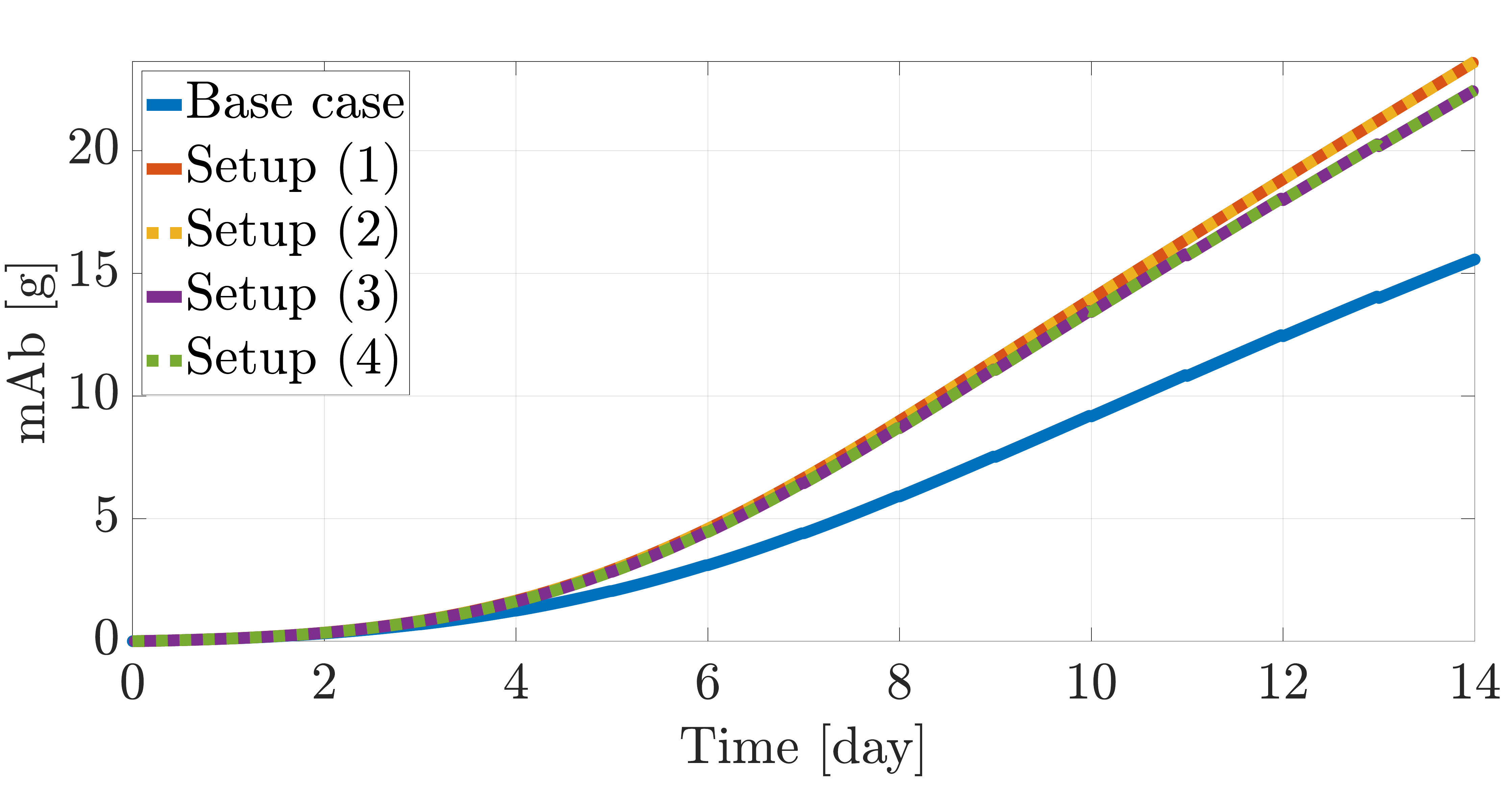}
    \caption{ Final mAb production for the base case simulation and the four optimal simulations. }
    \label{fig:mAbProduction}
\end{figure}

\subsection{ Optimization }
We solve the OCP, (\ref{eq:OPC}), for the full horizon of $T_h = 14$ day with $T_s = 30$ min. As such, we get the discrete horizon, $N = 672$. We consider four different optimization setups. We either apply no sampling or the sampling strategy from the base case. Additionally, we test the effect of enforcing almost full glucose utilization in the end of the fermentation. Table \ref{tab:optimizationSetups} presents the four optimization setups including the additional OCP constraints in each setup. We apply a direct multiple shooting discretization approach of the OCP and solve the resulting nonlinear programming (NLP) in CasADi \citep{andersson:2019a}. 

Fig. \ref{fig:optimalSolution} compares the base case simulation to the four optimal simulations. We observe that optimal cell growth is achieved by maintaining a constant glucose concentration of around $c_G = 7$ g/L to avoid glucose inhibition. If sampling is required, the optimal feeding trajectories compensates for the loss of medium to maintain the optimal glucose concentration. Once the product concentration fully inhibits the cell growth, the temperature is decreased to increase product formation. 




Fig. \ref{fig:mAbProduction} presents the time evolution of mAb in the bioreactor for the base case simulation and the optimal simulations. Table \ref{tab:mAbProduction} presents the total mAb production in the end of the $14$ day fermentation for all five simulations. We observe that all optimal simulations increase the mAb production compared to the base case. In particular, the optimal simulations produces $44$\% more mAb with sampling and $52$\% without sampling. We notice that forcing full glucose utilization does not affect the final production of mAb. This shows that there are multiple optimal solutions all resulting in the same mAb production \citep{ryde:2021a}. As such, optimization setup (2) and (4) are preferred under the assumption that glucose has a cost. This cost can be included in the optimization problem resulting in an economic OCP. 



%% file: tex/conclusion.tex
\section{Conclusion}
\label{sec:conclusion}

The paper presented a dynamic optimization numerical case study for mAb production. We applied a general modeling methodology to present an existing fermentation model for mAb production conducted in a continuous perfusion reactor. We expressed an optimization problem in terms of an OCP for maximization of mAb production in the end of the fermentation. Our results showed that optimal operation of the continuous perfusion reactor improves the mAb production by up to $52$\% compared to the base case. Additionally, our results showed that there exist multiple optimal solutions producing the same amount of mAb. Therefore, a full glucose utilization constraint in the OCP is preferred to reduce glucose loss.

%

